\documentclass[11pt]{elsarticle}
\usepackage{graphicx}
\usepackage{epsfig}
\usepackage{amssymb,amsmath}
\usepackage{textcomp}
\usepackage{wasysym}
\usepackage{psfrag}
\usepackage{a4wide}
\usepackage{wrapfig}
\usepackage{color}
\usepackage{bm}
\usepackage{array}
\usepackage{comment}

\journal{}
\begin{document}

\newcommand{\kv}[0]{\mathbf{k}}
\newcommand{\Av}[0]{\mathbf{A}}
\newcommand{\mv}[0]{\mathbf{m}}
\newcommand{\ev}[0]{\mathbf{e}}
\newcommand{\al}[0]{\mathbf{a_{1}}}
\newcommand{\as}[0]{\mathbf{a_{2}}}
\newcommand{\Mv}[0]{\mathbf{M}}
\newcommand{\Bv}[0]{\mathbf{B}}
\newcommand{\Iv}[0]{\mathbf{I}}
\newcommand{\dkx}[0]{\delta k_{x}}
\newcommand{\dky}[0]{\delta k_{y}}
\newcommand{\dk}[0]{\delta k}
\newcommand{\wv}[0]{\mathbf{w}}
\newcommand{\yv}[0]{\mathbf{y}}
\newcommand{\Rr}[0]{\Rv_{\rv}}
\newcommand{\xv}[0]{\mathbf{x}}
\newcommand{\vv}[0]{\mathbf{v}}
\newcommand{\Vv}[0]{\mathbf{V}}
\newcommand{\Sv}[0]{\mathbf{S}}
\newcommand{\rv}[0]{\mathbf{r}}
\newcommand{\Pv}[0]{\mathbf{P}}
\newcommand{\lv}[0]{\bm{\ell}}

\setlength{\jot}{2mm}

\newcommand{\jav}[1]{{\color{red}#1}}

\begin{frontmatter}

\title{Exactly solvable quadratic differential equation systems through generalized inversion}

\author{\'Ad\'am B\'acsi}
\ead{bacsi.adam@sze.hu}
\address{Department of Mathematics and Computational Sciences, Sz\'echenyi Istv\'an University \\ H-9026 Gy\H or, Egyetem t\'er 1, Hungary}
\address{MTA-BME Lend\"ulet Topology and Correlation Research Group,
Budapest University of Technology and Economics, \\ 1521 Budapest, Hungary}
\author{Albert Tiham\'er Kocsis}
\address{Department of Mathematics and Computational Sciences, Sz\'echenyi Istv\'an University \\ H-9026 Gy\H or, Egyetem t\'er 1, Hungary}

\date{\today}

\begin{abstract}
We study the autonomous systems of quadratic differential equations of the form $\dot{x}_i(t)=\xv(t)^T \Av_i \xv(t) + \vv_i^T \xv(t)$ with $\xv(t) = (x_1(t),x_2(t),\dots,x_i(t),\dots)$ which, in general, cannot be solved exactly. In the present paper, we present a subclass of analytically solvable quadratic systems, whose solution is realized through a multi-dimensional generalization of the inversion which transforms a quadratic system into a linear system. We provide a constructive algorithm which, on one hand, decides whether the system of differential equations is analytically solvable with the inversion transformation and, on the other hand, provides the solution. The presented results apply for arbitrary, finite number of variables.
\end{abstract}



\end{frontmatter}

\section{Introduction}
Linear, ordinary, autonomous systems of differential equations are celebrated for their analytical integrability which makes them excellent tools in almost all fields of science. This solvability breaks down when entering the realm of non-linear ordinary differential equation systems where no general solution method exists. In most cases, non-linear systems are handled with sophisticated numerical methods such as Runge-Kutta schemes\cite{iserles}.

Despite of the great accuracy of numerical methods, it is also desirable to find analytical techniques to solve non-linear systems. Beside the inherent beauty of exact solutions, they can also play the role of benchmark for numerical methods. In this work, we focus on first-order, quadratic, autonomous differential equation systems (QDEs) which emerge in various research fields including models of population dynamics\cite{champagnet_LV2010}, fluid dynamics\cite{coppel1991}, control systems \cite{control_elliott_2009} and even quantum dynamics \cite{bacsi2020}. Classification of two-variable quadratic systems have been studied extensively \cite{tsutomu,schlomiuk2005,CHAVARRIGA2006,reyn}.

There exists no general method to find analytical solution of QDEs. Except for the special case, when the number of variables is just one. By denoting this variable with $x(t)$, the quadratic differential equation is written as
\begin{gather}
\dot{x}=ax^2 + vx
\label{eq:ex1}
\end{gather}
with $a\neq 0$ and $v$ being constants. An additional constant term on the right-handside could be eliminated by shifting $x(t)$.
By introducing $y(t) = x(t)^{-1}$, the equation is transformed to the linear equation of
\begin{gather}
\dot{y}=-a-vy
\end{gather}
and, hence, the quadratic differential equation \eqref{eq:ex1} becomes exactly solvable. 

In this paper, we investigate a class of QDEs with more than one variable which can be solved exactly by transforming the QDE to a linear differential equation system (LDE). The transformation is realized through a multi-variable extension of $y(t) = x(t)^{-1}$, similar to the generalized spherical inversions presented in Refs. \cite{ramirez2014,ramirez2016}.



\section{Generalized inversion transformation}
Let us consider some differentiable functions $x_1(t)$, $x_2(t)$, \dots, $x_n(t)$ whose dynamics is governed by the QDE
\begin{gather}
\dot{x}_i = \xv^{T} \Av_i \xv + \vv_i^{T}\xv
\label{eq:diff}
\end{gather}
where 
$\Av_i$ are $n\times n$ matrices which are chosen to be symmetric and $\vv_i^T$ are row vectors of $n$ components. The entries of $\Av_i$ and $\vv_i$ are independent from $t$ and are assumed to be real without loss of generality\footnote{In case of complex entries, the differential equations can be split up to their real and imaginary parts.}. Our goal is to determine if the variables $x_i(t)$ can be transformed with a multivariable variant of the inversion in such a way that the QDE \eqref{eq:diff} is transformed to an LDE. The multivariable inversion is defined as
\begin{gather}
y_i(t) = \frac{x_i(t)}{\xv(t)^T\mathbf{B}\xv(t)}
\label{eq:btr}
\end{gather}
with some symmetric matrix $\Bv$ which does not depend on time. Note that we use the same $\Bv$ matrix for all $i$. This ensures that the transformation is easily inverted by
\begin{gather}
x_i(t) = \frac{y_i(t)}{\yv(t)^T\mathbf{B}\yv(t)}\,.
\label{eq:btrinverse}
\end{gather}
The transformations \eqref{eq:btr} and \eqref{eq:btrinverse} will henceforth be referred to as $B$-transformation and inverse $B$-transformation, respectively. 
The variables $y_i(t)$ are expected to obey the linear differential equations
\begin{gather}
\dot{y_i} = \mv_i^T\yv + w_i\,.
\label{eq:yeq}
\end{gather}
with time-independent row vectors $\mv_i^T$ and the scalars $w_i$. Since Eq. \eqref{eq:yeq} is exactly solvable for arbitrary initial conditions, the QDE can also be solved by applying the inverse $B$-transformation \eqref{eq:btrinverse}.
The goal of the present paper is to provide a procedure based on which one can decide whether a given QDE, defined through $\Av_i$ and $\vv_i$, can be $B$-transformed to an LDE. 



To start our study, let us calculate the time derivative of Eq. \eqref{eq:btrinverse}
\begin{gather}
\dot{x}_i=\frac{\dot{y}_i}{\yv^T\mathbf{B}\yv}  - \frac{y_i}{\left(\yv^T\mathbf{B}\yv\right)^2}\big(\dot{\yv}^T\Bv\yv + \yv^T\Bv\dot{\yv}\big)\,.
\end{gather}
Substituting the LDE differential equations \eqref{eq:yeq} and then transforming all $y_i$ variables back to $x_i$ by means of Eq. \eqref{eq:btr}, we obtain 
\begin{gather}
\dot{x}_i = \mv_i^{T}\xv + w_i\left(\xv^T\Bv\xv\right)- \xv^T\Bv\wv x_i -x_i\wv^T\Bv\xv - x_i\frac{\xv^T\left(\Mv^T\Bv +\Bv \Mv\right)\xv}{\xv^T\mathbf{B}\xv} 
\label{eq:dotxi}
\end{gather}
with the matrix $\Mv$ built from the row vectors $\mv_i^T$ and the vector $\wv$ consisting of the constants $w_i$.

By comparing Eq. \eqref{eq:dotxi} with the original QDE \eqref{eq:diff}, one can note that Eq. \eqref{eq:dotxi} is not necessarily a quadratic differential equation. The last term has an $\xv$-dependent denominator and, hence, is neither quadratic nor linear in general. This term scales with $\xv$ as a linear term but becomes an actual linear term only if
\begin{gather}
\Mv^T\Bv + \Bv\Mv = -\lambda\Bv
\label{eq:Mcond}
\end{gather}
holds with some scalar $-\lambda$. The minus sign is chosen for later convenience. The condition \eqref{eq:Mcond} means that $\Bv$ must be the eigenmatrix\cite{kronecker2009} of $\Mv$ with the eigenvalue of $-\lambda$. In this case, the differential equations \eqref{eq:dotxi} become
\begin{gather}
\dot{x}_i = \left(\mv_i^{T} +\lambda \ev_i^T\right)\xv + \xv^T\left( w_i \Bv-\Bv\left(\wv\circ\ev_i^T\right) - \left(\ev_i\circ\wv^T\right)\Bv \right) \xv
\label{eq:xeqB}
\end{gather}
where $\circ$ denotes dyadic product and $\ev_i^T$ is the unit row vector with $1$ in the $i$th entry and zero otherwise.
The equation is already a quadratic differential equation containing quadratic and linear terms in $\xv$. 

Comparing the linear terms of Eqs. \eqref{eq:xeqB} and \eqref{eq:diff}, we find
\begin{gather}
\Vv = \Mv + \lambda\Iv
\label{eq:vm}
\end{gather}
where $\Vv$ is constructed from the row vectors $\vv_i^T$ 
and $\Iv$ is the $n\times n$ identity matrix.
Substituting \eqref{eq:vm} into the condition \eqref{eq:Mcond}, we obtain
\begin{gather}
\Vv^T\Bv + \Bv\Vv = \lambda\Bv
\label{eq:Vcond}
\end{gather}
indicating that the matrix $\Bv$ must be a symmetric eigenmatrix of $\Vv$ as well. 

Beside the linear terms, the quadratic terms of Eq. \eqref{eq:xeqB} and \eqref{eq:diff} must also equal. Since both $\Av_i$ and the kernel in the quadratic terms in Eq. \eqref{eq:xeqB} are symmetric, they must equal
\begin{gather}
\Av_i = w_i\Bv - \Bv\left(\wv\circ\ev_i^T\right) - \left(\ev_i\circ\wv^T\right)\Bv
\label{eq:Acond}
\end{gather}
for all $i$.
Note that the second and third terms are matrices full of zeros except for the $i$th row and $i$th column, respectively. Therefore, omitting the $i$th row and column in the equation leads to the proportionality condition
\begin{gather}
\tilde{\Av}_i^i = w_i \tilde{\Bv}^i
\label{eq:prop}
\end{gather}
where $\tilde{\Av}_i^i$ ($\tilde{\Bv}^i$) is the $(n-1)\times(n-1)$ matrix which is obtained from $\Av_i$ ($\Bv$) by skipping the $i$th row and column. For the $i$th column of Eq. \eqref{eq:Acond}, we have
\begin{gather}
\mathbf{a}_i^i \equiv \Av_i\ev_i = w_i \Bv \ev_i -\Bv\wv - \left(\ev_i^T\Bv\wv\right)\ev_i\,.
\label{eq:ai}
\end{gather}

To summarize, the QDE as given in Eq. \eqref{eq:diff} can be transformed with a $B$-transformation to an LDE if we find a symmetric matrix $\Bv$ fulfilling Eq. \eqref{eq:Vcond} with some constant $\lambda$ and, at the same time, obeying Eqs. \eqref{eq:prop} and \eqref{eq:ai} with some constants $w_i$.
The linear terms of the LDE can be calculated by expressing $\Vv$ from Eq. \eqref{eq:vm}.

\section{Algorithm}
\label{sec:algo}
Let us provide a detailed description of the procedure based on the findings of the previous section.
The starting point of the algorithm is Eq. \eqref{eq:diff} determined by the matrices $\Av_i$ and vectors $\vv_i$.

\begin{itemize}
\item[1.] As a first step, one has to build up the matrix $\Vv$ from the row vectors $\vv_i^T$ and find all symmetric eigenmatrices based on $\Vv^T \Bv + \Bv\Vv = \lambda\Bv$.
Note that the eigenmatrices are definite only up to an overall multiplication factor similarly to the case of usual eigenvectors.

To obtain the eigenmatrices, one may calculate the vector-eigenvalues $s_j$ of $\Vv^T$ and the corresponding right-handside eigenvectors $\rv_j$, i.e., $ \Vv^T\mathbf{r}_{j} = s_j\mathbf{r}_{j}$. The sum of the geometric multiplicity of all eigenvalues is denoted by $N$. The matrix-eigenvalues of $\Vv$ are given by $s_1 + s_1$, $s_1+s_2$, \dots, $s_1 + s_n$, $s_2+s_2$, $s_2 + s_3$, \dots ,$s_N+s_N$ and the eigenmatrix corresponding to $s_j + s_m$ is $\mathbf{P}_{jm}=\left(\mathbf{r}_{j}\circ \mathbf{r}_{m}^T + \mathbf{r}_{m}\circ \mathbf{r}_{j}^T\right)/2$. Note that the number of symmetric eigenmatrices is $N(N+1)/2$ which is utmost $n(n+1)/2$ when $n$ eigenvectors are found.

In case of degenerate matrix eigenvalues $\lambda$, all linear combinations within the subspace are potential $\Bv$ matrices. For an example, see Sec. \ref{sec:3d}.

\item[2.] 
For each potential eigenmatrix $\Bv$, the proportionality condition \eqref{eq:prop} must be checked for all $i$ with some constants $w_i$,
which may also be zero.

If Eq. \eqref{eq:prop} cannot be fulfilled with any potential $\Bv$ matrix, then the QDE is not solvable by $B$-transforming to an LDE.
If Eq. \eqref{eq:prop} is satisfied with a $\Bv$ matrix and some $w_i$ constants, one may proceed to Step 3.


\item[3.] For the potential pairs of $\Bv$ and $\wv$ surviving Step 2, one has to check Eq. \eqref{eq:ai} for all $i$.
If it holds true for all $i$, then the QDE \eqref{eq:diff} can be transformed to an LDE of the form of Eq. \eqref{eq:yeq} where the constants $w_i$ are the coefficients of proportionality found in Eq. \eqref{eq:prop} in Step 2 and the vectors $\mv_i^T$ are the row vectors of $\Mv = \Vv - \lambda\Iv$
with $\lambda$ the eigenvalue corresponding to the eigenmatrix $\Bv$.

\item[4.] The solution of the QDE can be obtained by solving first the linear equation \eqref{eq:yeq} as 
\begin{gather}
\yv(t) = e^{\Mv t}\yv_0 + \yv_p(t)
\end{gather}
where $\yv_0 = \xv_0/(\xv_0^T\Bv\xv_0)$ is the initial condition with $\xv_0 = \xv(t=0)$ and $\yv_p(t) = \left(e^{\Mv t}-\Iv\right)\Mv^{-1}\wv$ if $\Mv$ is invertible. If $\Mv$ is not invertible, the function $\yv_p(t)$ has linear time-dependence in the nullspace of $\Mv$. Using the inverse $B$ transformation, the solution of the QDE is obtained as
\begin{gather}
\xv(t)=\frac{e^{\Mv t}\xv_0 + \left(\xv_0^T\Bv \xv_0\right) \yv_p(t)}
{e^{-\lambda t} + 2\yv_p(t)^T\Bv e^{\Mv t}\xv_0 + \left(\xv_0^T\Bv \xv_0\right)\left(\yv_p(t)^T\Bv \yv_p(t)\right)} \,.
\label{eq:sol}
\end{gather}
\end{itemize} 

Before presenting examples, let us make a remark on the region where the $B$-transformation is undefined, i.e., where $\xv^T\Bv\xv = 0$. This region can be any quadratic hypersurface depending on the specific structure of the $\Bv$ matrix. The region also involves the nullspace of $\Bv$, i.e., the points where $\Bv\xv=0$. In Sec. \ref{sec:2d}, the hypersurface contains two straight lines while in Sec. \ref{sec:3d}, it consists of two cone-shaped surfaces. This raises the question how the system behaves if the initial condition $\xv_0$ is an element of the region, i.e. $\xv_0^T\Bv\xv_0=0$. 

First, one can prove that if the initial condition of the QDE is on the hypersurface, the dynamics is constrained to the hypersurface  for the whole time-evolution. To justify this, we define $b(t) = \xv(t)^T\Bv\xv(t)$ whose dynamics is derived from the quadratic equations \eqref{eq:diff} as $\dot{b} = b\left(\lambda - 2\wv^T\Bv\xv\right)$. If the initial condition is on the hypersurface, i.e., $b(0)=0$, then the differential equation solves to $b(t)=0$. 

Second, it can also be shown that inside the hypersurface, the quadratic differential equation is analytically solvable and the solution is obtained simply by taking $\xv_0^T\Bv\xv_0=0$ in Eq. \eqref{eq:sol} which reads as
\begin{gather}
\xv(t)=\frac{e^{\Mv t}\xv_0}
{e^{-\lambda t} + 2\yv_p(t)^T\Bv e^{\Mv t}\xv_0}\,.
\label{eq:sol0}
\end{gather}
Note that \eqref{eq:sol0} does not obviously solve the quadratic differential equation because it was derived from Eq. \eqref{eq:sol} which was computed through the $B$ transformation. However, by substituting \eqref{eq:sol0} directly into \eqref{eq:diff} and taking advantage of $\xv_0^T\Bv\xv_0 = 0$ and the properties of $\Bv$, one can prove that the QDE is solved indeed by \eqref{eq:sol0}.
Hence, if a QDE is exactly solvable by a $B$-transformation, the integrability is preserved also in the region where the generalized inversion is undefined.


\section{Example, two-dimensional quadratic system}
\label{sec:2d}
In this section, we consider the quadratic equations of
\begin{align}
\dot{x}_1 & = -x_1^2 + 4x_1x_2 -x_1 +2x_2 \nonumber \\
\dot{x}_2 & = x_1^2 +2x_2^2 + x_1
\label{eq:example}
\end{align}
from which we read out
\begin{gather}
\Av_1 = \left[\begin{array}{cc} -1 & 2 \\ 2 & 0 \end{array}\right] \qquad
\Av_2 = \left[\begin{array}{cc} 1 & 0 \\ 0 & 2 \end{array}\right] \qquad
\vv_1^T = \left[\begin{array}{cc} -1 & 2 \end{array}\right] \qquad
\vv_2^T = \left[\begin{array}{cc} 1 & 0 \end{array}\right] \,.
\end{gather}
The first step is to obtain the eigenmatrices of $\Vv$. The right-handside eigenvectors of $\Vv^T$ are given by
\begin{gather}
\rv_1 = \left[\begin{array}{c} 2 \\ -2 \end{array}\right]\qquad
\rv_2 = \left[\begin{array}{c} 1 \\ 2 \end{array}\right]
\end{gather}
with the eigenvalues of $s_1 = -2$ and $s_2 = 1$. The potential $\Bv$ matrices are given by
\begin{align*}
\Pv_{11} &= \rv_1 \circ \rv_1^T = \left[\begin{array}{cc} 4 & -4 \\ -4 & 4 \end{array}\right] \qquad &\lambda_{11} = -4 \nonumber \\
\Pv_{22} &= \rv_2 \circ \rv_2^T = \left[\begin{array}{cc} 1 & 2 \\ 2 & 4 \end{array}\right] \qquad &\lambda_{22} = 2 \nonumber \\
\Pv_{12} &= \frac{1}{2}\left( \rv_1 \circ \rv_2^T + \rv_2 \circ \rv_1^T\right)  = \left[\begin{array}{cc} 2 & 1 \\ 1 & -4 \end{array}\right] \qquad &\lambda_{12} = -1
\end{align*}
which are non-degenerate eigenmatrices.

The second step is to check the proportionality condition of Eq. \eqref{eq:prop} and determine the coefficients $w_i$.

\renewcommand\arraystretch{2}


\begin{center}
\begin{tabular}{ | m{3cm} | m{4cm}| m{4cm} |} 
 \hline
  Proportionality check, Eq. \eqref{eq:prop} & $\tilde{A}_1^1 = [0]$ & $\tilde{A}_2^2 = [1]$ \\ 
  \hline
  $\Bv = \Pv_{11}$ & $\tilde{\Bv}^1=[1] \rightarrow w_1 = 0$ & $\tilde{\Bv}^2=[4] \rightarrow w_2 = \frac{1}{4}$ \\ 
  \hline
  $\Bv = \Pv_{22}$ & $\tilde{\Bv}^1=[4] \rightarrow w_1 = 0$ & $\tilde{\Bv}^2=[1] \rightarrow w_2 = 1$ \\ 
  \hline
  $\Bv = \Pv_{12}$ & $\tilde{\Bv}^1=[-4] \rightarrow w_1 = 0$ & $\tilde{\Bv}^2=[2] \rightarrow w_2 = \frac{1}{2}$ \\
  \hline
\end{tabular}
\end{center}

The table shows that all eigenmatrices obey the proportionality condition \eqref{eq:prop}. Note that for a two-dimensional QDE, the proportionality check always simplifies to comparison of $1\times 1$ matrices, which is in most cases trivially fulfilled. For higher dimensions, however, \eqref{eq:prop} might mean a much stricter condition. For an example, see Sec. \ref{sec:3d}.

In the example, we continue with the third step by checking Eq. \eqref{eq:ai} for each $i$. In the table below, $rhs$ stands for the right-handside of Eq. \eqref{eq:ai} evaluated with $\Bv$ and $\wv$.


\begin{center}
\begin{tabular}{ | m{3cm} | m{4cm}| m{4cm} |} 
  \hline
  $\mathbf{a}_i^i$ check, Eq.\,\eqref{eq:ai} & \renewcommand\arraystretch{1}
  $\mathbf{a}_1^1 = \left[\begin{array}{c} -1 \\ 2 \end{array}\right] $ & \renewcommand\arraystretch{1}
  $\mathbf{a}_2^2 = \left[\begin{array}{c} 0 \\ 2 \end{array}\right] $  \\ 
  \hline
 $\Bv = \Pv_{11}$ &
\renewcommand\arraystretch{1}
 $rhs = \left[\begin{array}{c} 2 \\ -1 \end{array}\right]\quad {\color{red} \times}$ & 
 \renewcommand\arraystretch{1}
$rhs = \left[\begin{array}{c} 0 \\ -1 \end{array}\right]\quad {\color{red} \times}$
\renewcommand\arraystretch{2} \\ 
  \hline
 $\Bv = \Pv_{22}$ &
\renewcommand\arraystretch{1}
$rhs = \left[\begin{array}{c} -4 \\ -4 \end{array}\right]\quad {\color{red} \times}$ & 
\renewcommand\arraystretch{1}
$rhs = \left[\begin{array}{c} 0 \\ -4 \end{array}\right]\quad {\color{red} \times}$ 
\renewcommand\arraystretch{2} \\ 
  \hline
  $\Bv = \Pv_{12}$  &  \renewcommand\arraystretch{1}
$rhs = \left[\begin{array}{c} -1 \\ 2 \end{array}\right]\quad {\color{green} \checkmark}$ & \renewcommand\arraystretch{1}
$rhs = \left[\begin{array}{c} 0 \\ 2 \end{array}\right]\quad {\color{green} \checkmark}$ 
\renewcommand\arraystretch{2}\\
  \hline
\end{tabular}
\end{center}

The investigation shows that $\mathbf{a}_1^1$ and $\mathbf{a}_2^2$ are reproduced only by $\Pv_{12}$. 

The result of the algorithm is that the QDE can be transformed to an LDE by $\Bv = \Pv_{12}$. The linear system is obtained as
\renewcommand\arraystretch{1}
\begin{gather}
\Mv = \Vv - \lambda_{12}\Iv = \left[\begin{array}{cc} 0 & 2 \\ 1 & 1 \end{array}\right] \qquad \wv = \left[\begin{array}{c} 0 \\ 1/2 \end{array}\right]\qquad\qquad
\begin{array}{rcl}
\dot{y}_1 &=& 2y_2 \\
\dot{y}_2 &=& y_1 + y_2 + \frac{1}{2}
\end{array}
\end{gather}
which can be solved analytically for arbitrary initial conditions. 

\begin{figure}[h]
\vspace*{-1cm}
\centering
\includegraphics[width=10cm]{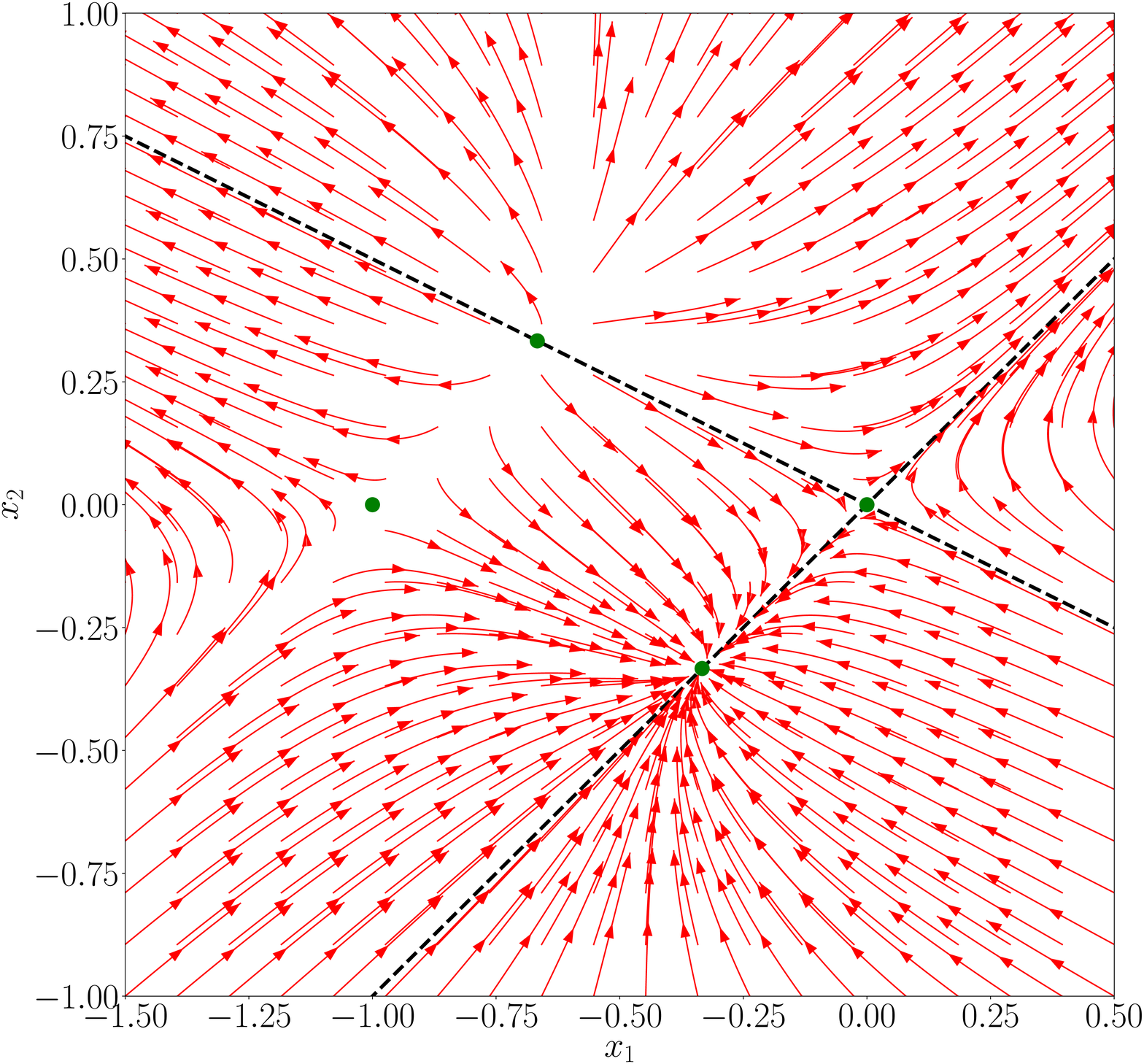}
\vspace*{-0.8cm}
\caption{Phase portrait of the quadratic system, Eq. \eqref{eq:example}. The green dots indicate the fix points of the system. The  dashed lines comprise the region where the $B$ transformation is undefined.}
\label{fig:phase}
\end{figure}

The analytical solution of the quadratic system can be given based on Eq. \eqref{eq:sol}. 
The phase portrait is shown in Fig. \ref{fig:phase}. 
The two dashed lines comprise the region where the $B$ transformation is undefined, i.e., where $\xv^T \Bv \xv = (\xv^T\rv_1)\cdot(\xv^T\rv_2)=0$. 

\newpage
\section{Three-dimensional example}
\label{sec:3d}
In this section, the algorithm is demonstrated in an example with three variables $x_1(t)$, $x_2(t)$ and $x_3(t)$. The quadratic differential equations are given by
\begin{gather}
\begin{array}{ccll} 
\dot{x}_1 & = & x_1^2 +7x_2^2 - 4x_1x_2 & + 5x_1  \\ & & & \\
\dot{x}_2  & =  & -2x_2^2+x_1x_2 + 2x_1 x_3 - 7 x_2x_3 & + 2x_2  \\ & & & \\
\dot{x}_3 & = & -x_2^2-7x_3^2 - 4x_2x_3 &  - x_3\end{array}
\label{eq:example3d}
\end{gather}
from which we can read out
\begin{gather}
\Av_1 =\left[\begin{array}{ccc} 1 & -2 & 0 \\ -2 & 7 & 0 \\ 0 & 0 & 0\end{array}\right] \qquad
\Av_2 =\left[\begin{array}{ccc} 0 & \frac{1}{2} & 1 \\ \frac{1}{2} & -2 & -\frac{7}{2} \\ 1 & -\frac{7}{2} & 0\end{array}\right] \qquad
\Av_3 =\left[\begin{array}{ccc} 0 & 0 & 0 \\ 0 & -1 & -2 \\ 0 & -2 & -7\end{array}\right]
\end{gather}
and
\begin{gather}
\Vv =\left[\begin{array}{ccc} 5 & 0 & 0 \\ 0 & 2 & 0 \\ 0 & 0 & -1\end{array}\right]\,.
\end{gather}
The eigenvectors of $\Vv^T$ are simply obtained as $\rv_i = \ev_i$ with the eigenvalues of $s_1 = 5$, $s_2=2$ and $s_3=-1$.
Hence, the eigenmatrices of $\Vv$ matrices are as follows.
\begin{align*}
\Pv_{11} &= \left[\begin{array}{ccc} 1 & 0 & 0 \\ 0 & 0 & 0 \\ 0 & 0 & 0 \end{array}\right] \qquad &\lambda_{11} = 10 \nonumber \\
\Pv_{22} &= \left[\begin{array}{ccc} 0 & 0 & 0 \\ 0 & 1 & 0 \\ 0 & 0 & 0 \end{array}\right] \qquad &\lambda_{22} = 4 \nonumber \\
\Pv_{33} &= \left[\begin{array}{ccc} 0 & 0 & 0 \\ 0 & 0 & 0 \\ 0 & 0 & 1 \end{array}\right] \qquad &\lambda_{33} = -2 \nonumber \\
\Pv_{12} &= \left[\begin{array}{ccc} 0 & \frac{1}{2} & 0 \\ \frac{1}{2} & 0 & 0 \\ 0 & 0 & 0 \end{array}\right] \qquad &\lambda_{12} = 7 \nonumber \\
\Pv_{13} &= \left[\begin{array}{ccc} 0 & 0 & \frac{1}{2} \\ 0 & 0 & 0 \\ \frac{1}{2} & 0 & 0 \end{array}\right] \qquad &\lambda_{13} = 4 \nonumber \\
\Pv_{23} &= \left[\begin{array}{ccc} 0 & 0 & 0 \\ 0 & 0 & \frac{1}{2} \\ 0 & \frac{1}{2} & 0 \end{array}\right] \qquad &\lambda_{23} = 1 
\end{align*}
Note that the eigenmatrices $\Pv_{22}$ and $\Pv_{13}$ are degenerate. Therefore, their linear combination
\begin{gather}
\Pv_4 = \Pv_{22} + b\Pv_{13} = \left[\begin{array}{ccc} 0 & 0 & \frac{b}{2} \\ 0 & 1 & 0 \\ \frac{b}{2} & 0 & 0 \end{array}\right]
\end{gather}
is also a potential $\Bv$ matrix with an arbitrary value of $b$.

We continue the algorithm with Step 2, the proportionality check. 

\renewcommand\arraystretch{2}
\vspace{-0.5cm}
\begin{center}
\hspace*{-1cm}\begin{tabular}{ | m{2.2cm} | m{4.6cm}| m{4.6cm} | m{4.6cm} |}
 \hline
  Prop. check, Eq. \eqref{eq:prop} & 
  \renewcommand\arraystretch{1} $\tilde{A}_1^1 = \left[\begin{array}{cc} 7 & 0 \\ 0 & 0 \end{array} \right]$ & 
  \renewcommand\arraystretch{1} $\tilde{A}_2^2 = \left[\begin{array}{cc} 0 & 1 \\ 1 & 0 \end{array} \right]$ & 
  \renewcommand\arraystretch{1} $\tilde{A}_3^3 = \left[\begin{array}{cc} 0 & 0 \\ 0 & -1 \end{array} \right]$ \renewcommand\arraystretch{2}\\ 
  \hline
  $\Bv = \Pv_{11}$ & 
  \renewcommand\arraystretch{1} $\tilde{\Bv}^1=\left[\begin{array}{cc} 0 & 0 \\ 0 & 0 \end{array} \right] \rightarrow w_1 = 0$ &
  \renewcommand\arraystretch{1} $\tilde{\Bv}^2=\left[\begin{array}{cc} 1 & 0 \\ 0 & 0 \end{array} \right] \qquad  {\color{red} \times}$ & \qquad -- \renewcommand\arraystretch{2}\\
  \hline
  $\Bv = \Pv_{33}$ & 
  \renewcommand\arraystretch{1} $\tilde{\Bv}^1=\left[\begin{array}{cc} 0 & 0 \\ 0 & 1 \end{array} \right] \qquad  {\color{red} \times}$ &
  \qquad -- & \qquad -- \renewcommand\arraystretch{2}\\
  \hline
  $\Bv = \Pv_{12}$ & 
  \renewcommand\arraystretch{1} $\tilde{\Bv}^1=\left[\begin{array}{cc} 0 & 0 \\ 0 & 0 \end{array} \right] \rightarrow w_1 = 0$ &
  \renewcommand\arraystretch{1} $\tilde{\Bv}^2=\left[\begin{array}{cc} 0 & 0 \\ 0 & 0 \end{array} \right] \rightarrow w_2 = 0$ & 
  \renewcommand\arraystretch{1} $\tilde{\Bv}^3=\left[\begin{array}{cc} 0 & \frac{1}{2} \\ \frac{1}{2} & 0 \end{array} \right] \qquad  {\color{red} \times}$ \renewcommand\arraystretch{4}\\
  \hline
  $\Bv = \Pv_{23}$ & 
  \renewcommand\arraystretch{1} $\tilde{\Bv}^1=\left[\begin{array}{cc} 0 & \frac{1}{2} \\ \frac{1}{2} & 0 \end{array} \right] \qquad  {\color{red} \times}$ &  \qquad -- &   \qquad -- \renewcommand\arraystretch{2}\\
  \hline
  $\Bv = \Pv_{4}$ & 
  \renewcommand\arraystretch{1} $\tilde{\Bv}^1=\left[\begin{array}{cc} 1 & 0 \\ 0 & 0 \end{array} \right] \rightarrow w_1 = 7$ &
  \renewcommand\arraystretch{1} $\tilde{\Bv}^2=\left[\begin{array}{cc} 0 & \frac{b}{2} \\ \frac{b}{2} & 0 \end{array} \right] \rightarrow w_2 = \frac{2}{b}$ & 
  \renewcommand\arraystretch{1} $\tilde{\Bv}^3=\left[\begin{array}{cc} 0 & 0 \\ 0 & 1 \end{array} \right] \rightarrow w_3 = -1$ \renewcommand\arraystretch{2}\\
  \hline
\end{tabular}
\end{center}
Note that only $\Pv_4$ obeys the proportionality check but no restriction on $b$ has been obtained.
Hence, in Step 3, we only investigate $\Pv_4$.
\renewcommand\arraystretch{2}
\vspace{-0.5cm}
\begin{center}
\hspace*{-1cm}\begin{tabular}{ | m{2.2cm} | m{4.6cm}| m{4.6cm} | m{4.6cm} |}
 \hline
  $\mathbf{a}_i^i$ check, Eq.\,\eqref{eq:ai} & 
  \renewcommand\arraystretch{1} $\renewcommand\arraystretch{1}\tilde{a}_1^1 = \left[\begin{array}{c}  1 \\ -2 \\ 0 \end{array} \right]$ & 
  \renewcommand\arraystretch{1} $\tilde{a}_2^2 = \left[\begin{array}{c} \frac{1}{2} \\ -2 \\ -\frac{7}{2} \end{array} \right]$ & 
  \renewcommand\arraystretch{1} $\tilde{a}_3^3 = \left[\begin{array}{c} 0 \\ -2 \\ -7 \end{array} \right]$ \renewcommand\arraystretch{2}\\ 
  \hline
  $\Bv = \Pv_{4}$ & \renewcommand\arraystretch{1}
  $ rhs = \left[\begin{array}{c} b \\ -\frac{2}{b} \\ 0\end{array} \right] {\color{green} \checkmark} \quad\mbox{if $b=1$}$ &
  $\renewcommand\arraystretch{1} rhs = \left[\begin{array}{c} \frac{b}{2} \\ -\frac{2}{b} \\ -\frac{7b}{2} \end{array} \right]  {\color{green} \checkmark} \, \mbox{if $b=1$}$ & 
  $\renewcommand\arraystretch{1} rhs = \left[\begin{array}{c} 0 \\ -\frac{2}{b} \\ -7b\end{array} \right] {\color{green} \checkmark}\, \mbox{if $b=1$}$ \\
    \hline
\end{tabular}
\end{center}
It has been found that the quadratic system transformed to a linear system with the matrix
\renewcommand\arraystretch{1}
\begin{gather}
\Bv = \left[\begin{array}{ccc} 0 & 0 & \frac{1}{2} \\ 0 & 1 & 0 \\ \frac{1}{2} & 0 & 0 \end{array}\right]
\end{gather}
and the resulting linear differential equations are given as
\begin{gather}
\begin{array}{rcl}
\dot{y}_1 &=& y_1 + 7 \\
\dot{y}_2 &=& -2y_2 + 2 \\
\dot{y}_3 &=& -5y_3 - 1
\end{array}
\end{gather}
which are exactly solvable for arbitrary initial conditions. The analytical solution can be given based on Eq. \eqref{eq:sol}.

\begin{figure}[h]
\centering
\includegraphics[width=7cm]{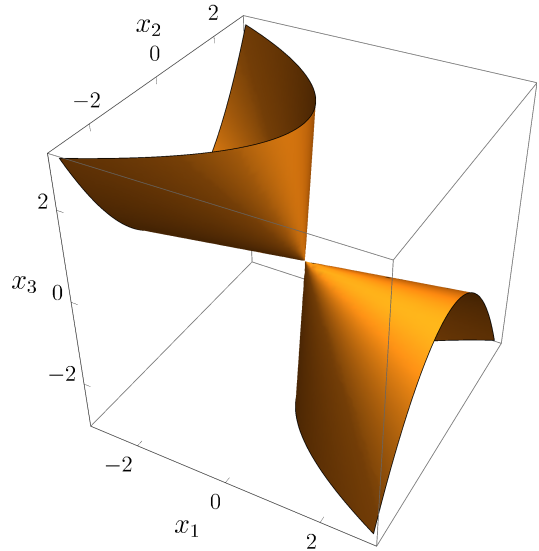}
\caption{The hypersurface determined by $\xv^T\Bv\xv=0$ consists of two cone-shaped surfaces. The hypersurface contains the two special directions determined by $\lv_1$ (blue) and $\lv_3$ (red).}
\label{fig:cones}
\end{figure}

Note that the $B$ transformation is undefined in the region where $\xv^T \Bv \xv = x_1x_3 + x_2^2=0$. This quadratic equation determines  two cone-shaped surfaces touching each other at the origin as shown in Fig. \ref{fig:cones}.

\section{Conclusion}
We studied the system of quadratic differential equations of the form of Eq. \eqref{eq:diff}. Although these differential equation systems cannot be solved in general, we have presented a procedure by means of which some quadratic systems can be solved analytically. The solution is performed through the $B$-transformation of $\yv = \xv/(\xv^T\Bv \xv)$ which is a multi-dimensional generalization of the inversion $y= x^{-1}$.

In Sec. \ref{sec:algo}, we have described the algorithm which allows one to decide whether a quadratic system can be transformed with a $B$-transformation to a linear differential equation system. If so, the algorithm also yields the linear system which is exactly solvable. The algorithm works for arbitrary number of variables.


\bibliographystyle{elsarticle-num} 
\bibliography{quadratic}

\newpage
\appendix

\end{document}